\documentclass[12pt]{amsart}
\usepackage{amsthm,amsmath,amssymb,amscd}
\newcommand{\thisdate}{\today}
{\newtheorem*{maintheorem}{Main Theorem}
   \newtheorem{theorem}{Theorem}[subsection]                     
   \newtheorem{proposition}[theorem]{Proposition}     
   \newtheorem{lemma}[theorem]{Lemma}
   \newtheorem*{claim}{Claim}

   \newtheorem{conjecture}[theorem]{Conjecture}
   
}
{\theoremstyle{definition}
   \newtheorem{definition}[theorem]{Definition}
   
}

\newcommand{\CC}{{\mathbb{C}}}
\newcommand{\QQ}{{\mathbb{Q}}}

\newcommand{\ZZ}{{\mathbb{Z}}}

\newcommand{\bA}{{\mathbf{A}}}

\newcommand{\bM}{{\mathbf{M}}}

\newcommand{\cN}{{\mathcal N}}
\newcommand{\cO}{{\mathcal O}}
\newcommand{\cP}{{\mathcal P}}
\newcommand{\cQ}{{\mathcal Q}}
\newcommand{\cR}{{\mathcal R}}

\newcommand{\cTheta}{{\tilde{{\varTheta}}}}

\newcommand{\fp}{{\mathfrak p}}
\newcommand{\fA}{{\boldsymbol{\mathfrak A}}}
\newcommand{\fXi}{{\boldsymbol{\Xi}}}
\newcommand{\fB}{{\boldsymbol{\mathfrak B}}}
\newcommand{\obA}{\overline{\mathbf A}}
\newcommand{\ofA}{\overline{\boldsymbol{\mathfrak A}}}
\newcommand{\obTheta}{{\overline{\boldsymbol{\Theta}}}}
\newcommand{\ofTheta}{{\overline{\tilde{\boldsymbol{\varTheta}}}}}
\newcommand{\Spec}{\operatorname{Spec}}
\newcommand{\Isom}{\operatorname{Isom}}
\newcommand{\Sing}{\operatorname{Sing}}

\newcommand{\Pic}{{\operatorname{\mathbf{Pic}}}}

\newcommand{\Aut}{{\operatorname{Aut}}}

\newcommand{\length}{\operatorname{length}}

\newcommand{\das}{\dashrightarrow}
\newcommand{\dar}{\downarrow}
\newcommand{\ocM}{\overline{{\mathcal M}}}
\newcommand{\obM}{\overline{{\mathbf M}}}

\newcommand{\setmin}{\,\protect%
\begin{picture}(8,10)\qbezier(1,5.5)(4,4.)(7,2.5)\end{picture}\,}

\begin{document}
\title[Uniformity of stably integral points (\thisdate)]{Uniformity of stably
   integral points on principally polarized abelian 
varieties of dimension $\leq 2$}
\author[D. Abramovich]{Dan Abramovich}
\thanks{D.A. Partially supported by NSF grant DMS-9700520 and by an Alfred
P. Sloan research fellowship}  
\address{Department of Mathematics\\ Boston University\\ 111 Cummington\\
Boston, MA 02215\\ USA} 
\email{abrmovic@math.bu.edu}
\author[K. Matsuki]{Kenji Matsuki}
\thanks{K.M. partially supported by NSA grant
		MDA904-96-1-0008.}
\address{Department of Mathematics \\ Purdue University \\  
                        1395 Mathematical Sciences Building \\  West Lafayette,
                        IN 47907-1395 \\ USA}
\email{kmatsuki@math.purdue.edu}
\date{\thisdate}

\begin{abstract} The purpose of this paper is to prove, assuming that the
conjecture of Lang and 
Vojta holds true, that there is a uniform bound on the number of stably
integral points in the complement of the theta divisor on a principally
polarized abelian surface defined over a number field.
 Most of our argument works
in arbitrary dimension and the restriction on the dimension $\leq 2$ is used
only at the last step, where we apply Pacelli's stronger
uniformity results for elliptic curves. 
\end{abstract}

\maketitle

{\bf Preliminary version, \thisdate.}

\addtocounter{section}{-1}
\section{Introduction}

\subsection{The conjecture of Lang and Vojta.}\label{Ssec:lgt} Let $X$ be a
variety 
over a field of characteristic 0. We say that $X$ is a variety of logarithmic
general 
type, if there exists a desingularization $\tilde{X} \to X$, and a projective
embedding $\tilde{X} \subset Y$ where $D = Y\setmin \tilde{X}$ is a divisor of
normal 
crossings, such that the invertible sheaf $\omega_Y(D)$ is big. We note first
that this peoperty is independent of the choices of $\tilde{X}$ and $Y$, and
that it is a {\em proper birational invariant}, namely, if $X'\to X$ is a
proper birational morphism (or an inverse of such) then $X$ is of logarithmic
general  type if and only if $X'$ is.

Now let $X$ be a variety of
logarithmic general type defined over a number 
field
$K$.  Let $S$ be a finite set of places in $K$ and let ${\mathcal
O}_{K,S}\subset K$ be the 
ring of $S$-integers. Fix a model ${\mathcal X}$  of $X$ over ${\mathcal
O}_{K,S}$.  It was conjectured by S. Lang and P. Vojta  (cf. \cite{lang:conj}, 
\cite{Vojta:higher-Mordell}) that the 
set of 
$S$-integral points on ${\mathcal X}$ is not Zariski dense in ${\mathcal X}$.
In 
case $X$ is projective, one may choose an arbitrary projective model ${\mathcal
X}$ and then ${\mathcal X}({\mathcal O}_{K,S})$ is identified with $X(K)$.  In
such a case, one often refers to this Lang-Vojta conjecture  as
{\em Lang's conjecture.} When $\dim X=1$, the conjectures of Lang and Vojta
reduce to Siegel's theorem and Mordell's conjecture (Faltings's theorem).

\subsection{The uniformity principle} \label{Sec:uniformity}\label{Subsec:principle} In
\cite{C-H-M}, L. Caparaso, 
J. Harris and B. Mazur show that Lang's conjecture  
implies a uniformity result for rational points on curves of genus $g \geq 2$
over a fixed number field, which extends Faltings's theorem
\cite{Faltings:Mordell}:

\begin{quote} 
{\em Suppose Lang's conjecture holds true. Then there exists a number 
$N(g,K)$ (depending only on the genus
$g$ and the number field $K$) such that the number of rational points $\# C(K)$
on a smooth projective curve $C$ of genus $g$ defined over a number field $K$
is 
uniformly bounded
$$\# C(K) < N(g,K).$$}
\end{quote}

 The basic
principle of \cite{C-H-M} may be 
summarized by the 
implication
%\begin{quote}
\bigskip

\parbox{1.5in}{Lang's conjecture in arbitrary dimension}
$\quad\Longrightarrow\quad $
\parbox{4.3in}{Uniform version of Lang's conjecture in a fixed dimension \\
(e.g. Uniform Mordell's
conjecture in dimension 1).}
%\end{quote}

\bigskip

Indeed, the results of \cite{Hassett:surfaces}, \cite{A:fibered} and
\cite{A-Voloch:fibered} show that 
the same principle holds in higher dimensions as well.

It is only natural then to seek to show the following analogous implication,
which may be considered a logarithmic generalization of the above:
\smallskip

%\begin{quote}
\parbox{2.in}{The Lang-Vojta conjecture in arbitrary dimension}
$\quad\Longrightarrow\quad $
\parbox{3.8in}{Uniform version of the Lang-Vojta conjecture\\ in a fixed
dimension.}

\subsection{The case of elliptic curves.} \label{Sec:intro-elliptic-case}
Let $K$ be a number field,  $S$ a finite number of places in $K$, and denote
by ${\mathcal O}_{K,S}$ the 
ring of $S$-integers. Let $E$ be an elliptic curve defined
over $K$, with origin $0$,  and let $P$ be a $K$-rational point of $E \setmin
\{0\}$, i.e., 
$P \in (E \setmin \{0\})(K)$. Fix  a model
${\mathcal E}$ of $E$ over ${\mathcal O}_{K,S}$, and denote by
$\overline{\{0\}}$ the ``zero section''.  We say that $P$ is
$S$-integral  if 
$P \in ({\mathcal E} \setmin
\overline{\{0\}})({\mathcal O}_{K,S})$.  Siegel's theorem, which may be
considered a 
logarithmic version of Mordell's conjecture, states that the number of
$S$-integral points on ${\mathcal E}\setmin \overline{\{0\}}$ is finite. 
We can view this as a case of the Lang-Vojta conjecture, 
 regarding $E \setmin \{0\}$ as a curve of logarithmic general type.
 Thus according to the 
principle in \ref{Subsec:principle}, one might naively expect that, assuming
the  Lang-Vojta conjecture, a uniform version
of Siegel's 
theorem in the following  form would hold:
\begin{quote} Could there exist a number $N(K,S)$ (depending only on the
number field and the finite 
number of places $S$ in $K$) such that for any elliptic curve defined over $K$
and any  model ${\mathcal E}$  over ${\mathcal O}_{K,S}$, the number of
$S$-integral points in the complement of the zero section is uniformly bounded:
$$\# ({\mathcal E} \setmin \overline{\{0\}})({\mathcal O}_{K,S}) <
N(K,S)\hbox{?}$$  
\end{quote}
{\em but the statement in this naive form fails to
hold.} Indeed,  take
an elliptic curve $E$ 
with an infinite number of $K$-rational points.  Let
$$y^2 = x^3 + Ax + B$$
be an affine equation for $E$ with $A, B \in {\mathcal O}_{K,S}$ (here the
origin of $E$ is the point at infinity). This equation 
gives an integral model ${\mathcal E}$ over ${\mathcal O}_{K,S}$.   Note
that for an arbitrary $n$-tuple of  $K$-rational 
points
$P_1, \cdot\cdot\cdot, P_n \in (E \setmin \{0\})(K)$, one can find $c \in
{\mathcal 
O}_{K,S}$ such that $c^2x(P_i), c^3y(P_i) \in {\mathcal O}_{K,S}$ where
$x(P_i)$ 
and $y(P_i)$ are $x$- and $y$-coordinates of the point $P_i$.  By changing
coordinates $x_1 = c^2x, y_1 = c^3y$, one obtains a new  model
${\mathcal 
E}'$ of
$E$ with a different defining equation over
${\mathcal O}_{K,S}$:
$$y_1^2 = x_1^3 + c^4Ax_1 + c^6B$$
where all the points $P_i$ are now $S$-integral points in ${\mathcal E}'$.
This 
example shows that {\em even for a fixed elliptic curve} defined over $K$ one
may 
have an arbitrarily large number of
$S$-integral points on varying  models over $\cO_{K,S}$, and hence the number
is not 
uniformly bounded.

We observe that this unboundedness is caused, as
demonstrated in the  example above, by allowing some coordinate
changes. Geometrically, these coordinate changes
correspond to some blowing up centered at the zero points in some fibers of
${\mathcal E}\to \Spec{\mathcal O}_{K,S}$, possibly followed
by some blowing down. From the Lang-Vojta point of view, such a procedure may
introduce a curve $F$ in a fiber with 
negative  
intersection with the logarithmic relative dualizing sheaf
$$\omega_{{\mathcal E}/\Spec{\mathcal O}_{K,S}}(\overline{\{0\}}) \cdot F \leq
0.$$  Such a component fails
to be ``hyperbolic'' and thus may ``leave space'' for more integral points.

In order to avoid such a situation one may wish to impose some positivity
condition on the 
 models one takes. This lead the first author to the notion of 
{\em stably} 
$S$-integral points:  a $K$-rational point
$P \in (E \setmin \{0\})(K)$ is called stably $S$-integral if for any finite
extension $L$ of $K$, with $T$ being the set of places above $S$, such that $E$
has a 
stable model ${\mathcal E}_{L,T}$ over ${\mathcal O}_{L,T}$, we have $P \in
({\mathcal 
E}_{L,T} \setmin \overline{\{0\}})({\mathcal O}_{L,T})$.  

%(Recall that ${\mathcal E}__{L,T}$ is
%a semi-stable model iff it is nonsingular, and any singular fiber of
%${\mathcal 
%E}_L \rightarrow Spec {\mathcal O}_{L,T}$ is reduced, its components crossing
%normally 
%(including the zero section) and the logarithmic relative dualizing
%sheaf $\omega_{{\mathcal E}_L/Spec {\mathcal O}_{L,T}}(\{0\})$ is relatively
%nef over 
%$Spec {\mathcal O}_{L,T}$.  A stable model is the (relative) log canonical
%model 
%of a semistable model over
%$Spec {\mathcal O}_{L,T}$ and it may be singular.)

Using this definition, the following assertion was  shown in
\cite{A:stintegral}: 
\begin{quote} 
Assume that the Lang-Vojta conjecture holds true. Then
there 
exists a number $N(K,S)$ such that for any elliptic curve defined over
$K$ the number of stably
$S$-integral points $E(K,S)^{stable}$ is uniformly
bounded:
$$\# E(K,S)^{stable} < N(K,S).$$
\end{quote}

\subsection{Abelian varieties.} The purpose of this paper is to extend the
result of \ref{Sec:intro-elliptic-case} to the higher 
dimensional case, according to the uniformity principle of
\ref{Subsec:principle}.  Let 
$A$ be a principally 
polarized abelian variety with  theta divisor $\Theta$, defined over a number
field $K$.  Let
$S$ be a finite set of places in $K$, and ${\mathcal O}_{K,S}$ the ring of
$S$-integers.  It is a theorem due to Faltings \cite{Faltings:abel} that if
${\mathcal A}\to {\mathcal O}_{K,S}$ is 
a model of
$A$ over ${\mathcal O}_{K,S}$, and $\overline{\varTheta}$ the closure of $\Theta$,
then 
$$\# 
({\mathcal A} \setmin \overline{\varTheta})({\mathcal O}_{K,S}) < \infty.$$  
As observed in \ref{Sec:intro-elliptic-case} 
for the case of dimension 1,  one cannot expect that the number of
$S$-integral points be  uniformly bounded
without imposing some positivity condition, or equivalently, without
restricting 
oneself to some notion of stably integral points. In Definition
\ref{Def:stably-integral} below we  define  stably
integral points as points which are integral 
on the complement of $\overline{\varTheta}$ on the stable model
of a principally polarized abelian
variety,  after taking a finite extension of $K$.  The existence of 
such stable models is provided by 
recent results of Alexeev and Nakamura
(see \cite{Alexeev:SQAV}, \cite{Alexeev:pairs}, \cite{Alexeev-Nakamura},
\cite{Alexeev:toric}), in 
which the 
moduli of principally polarized abelian varieties is compactified by the moduli
of stable quasi-abelian  pairs.

Here is our main arithmetic result on abelian surfaces:
\begin{maintheorem}
 Assume that the conjecture of Lang and Vojta holds true. Then there
exsits a number $N(K,S)$, such that for any principally polarized abelian
{\em surface} with a 
theta divisor $(A, \Theta)$, the number of the stably $S$-integral points of 
$A \setmin\Theta$ is uniformly bounded:
$$\# (A\setmin\Theta)(K,S)^{stable} < N(K,S).$$
\end{maintheorem}

We expect a similar result to hold for 
abelian varieties of any fixed dimension with a divisor of arbitrary fixed
polarization degree. 

\subsection{}
In Section \ref{Sec:review}, we review the proof of the  uniformity statement
on  rational points on
varieties of general type, as we will apply several methods which have
been used in that context. In Section \ref{Sec:Alex} we review the contruction
of complete moduli of stable quasi-abelian pairs by Alexeev and Nakamura. In
Section \ref{Sec:proof}, we   prove the Main 
Theorem.  There is one difficulty in the last inductive step
where we consider  {\em families of  subvarieties} of principally
polarized abelian varieties, especially families of {\em abelian} subvarieties,
which are not necessarily {\em principally} polarized.
 We can complete the argument in dimension
$\leq 2$ using Pacelli's stronger uniformity results for the elliptic curves,
leaving the general case of dimension
$> 2$ conjectural.

\subsection{} It is worth noting that a similar argument to the one we give for
principally polarized abelian varieties, often works for  pairs of logarithmic
general type  
$(X,D)$ (see \ref{Ssec:lgt-pair} below) defined over $K$, if a ``good'' moduli
for 
the log 
canonical models of such pairs exists.  For example, one can use such an
argument to show, assuming the conjecture of Lang and Vojta, that there is a
uniform bound on  stably integral points for ${\mathbb P}^1 \setmin \{n \text{\
points}\} 
\hskip.1in (n
\geq 3)$, using the moduli of stable $n$-pointed curves of genus 0.
However, at least  when $n = 3$, this uniformity statement  is nothing but the
classical result of Siegel about the finiteness 
of the number of solutions of an
$S$-unit equation (which holds regardless of the Lang-Vojta conjecture). This 
result of Siegel has been strengthened to a great 
extent in recent years.

\subsection{Logarithmic pairs}\label{Ssec:lgt-pair} In Section \ref{Ssec:lgt}
we defined what it 
means for a variety to be of logarithmic general type, in terms of a
good compactification $\tilde{X} \subset Y$ of a desingularization
$\tilde{X}\to X$.  
It is convenient to have a criterion which does not require choosing a
desingularization. One can approach that using ``singularities of pairs'', 
see \cite{KMM}.

\begin{quote} Let $\overline{X}$ be a projective variety, $D\subset
\overline{X}$ a reduced effective Weil divisor. 
 Assume
 \begin{enumerate}
\item  the pair $(\overline{X},D)$ has
log-canonical singularities;
\item 
$\omega_{\overline{X}}(D)$ is big, and 
\item the complement $\overline{X}\setmin D$ has canonical singularities. 
\end{enumerate}
Then $(\overline{X},D)$ has  logarithmic general type.
\end{quote}

Thus it is enough to check that $\overline{X}$ has a lot of logarithmic
differentials, and that its singularities are sufficiently mild.

We would like to draw the reader's attention to condition 3, which does not
follow from condition 1 because of possible exceptional divisors. Many authors
define the {\em pair} $(\overline{X},D)$ to be of logarithmic general type if
conditions 1 and 2 are satisfied - this is equivalent to the statement that the
 variety 
$\overline{X} \setmin (D \cup \Sing(X))$ is of logarithmic general type in our
terminology.  

\subsection{Aknowledgements} We are grateful to Professors Alexeev, Hassett,
Kawamata, Koll\'ar and Pacelli for  
invaluable suggestions on various parts of the paper.

%------------------------------------------------------------------------------
\section{Outline of the proof of uniformity for varieties of general type.}
\label{Sec:review}
%------------------------------------------------------------------------------
%

%In this section, we outline the proof of the uniformity after
%[CHM] and [H] for rational points on curves of genus
%$\geq 2$ and more in general varieties of general type of higher dimension
%under Lang's 
%conjecture, with some modifications toward \S 2 in mind.

\subsection{Correlation of points} We briefly recall the outline of the proof
of  uniformity of rational points on curves 
of genus $\geq 2$ in \cite{C-H-M}.  One of the main ideas of \cite{C-H-M} is to
observe that, assuming  Lang's conjecture holds true, the
set of
$K$-rational points on all smooth projective curves  of genus $g \geq 2$
defined over a number field $K$
is {\em correlated}, i.e., the collection of  $n$-tuples of such points
satisfies  a nontrivial algebraic relation, for suitable $n$.

Let $\pi:X \rightarrow B$ be a projective family of smooth irreducible curves
of genus $g \geq 2$ defined over a number field $K$.  We denote by
$$\pi_n:X_B^n = X \mathop{\times}\limits_B\ldots\mathop{\times}\limits_BX
\rightarrow B$$ the $n$-th fibered 
power of $X$ over $B$.  Denote 
by $\tau_n:X_B^n \rightarrow X_B^{n-1}$ the projection onto the first $n-1$
factors.  Given a point $b \in B$ we denote by $X_b$ the fiber
$\pi^{-1}(b)$.  Similarly, given a point $Q = (P_1, \cdot\cdot\cdot, P_{n-1})
\in X_B^{n-1}$ we denote by $X_Q \subset X_B^n$ the fiber
$\tau_n^{-1}(Q)$.  Note that if $Q\in X_B^{n-1}(K)$ and  $\pi_{n-1}(Q) = b$
then $X_Q \cong X_b$. 

Assume that we are given a subset ${\mathcal P} \subset X(K)$.  We denote by
${\mathcal P}_B^n \subset X_B^n$ the fibered power of ${\mathcal P}$ over $B$
(namely 
the union of the $n$-tuples of points in ${\mathcal P}$ whose images on $B$ are
the same), and by ${\mathcal 
P}_b$ the points of
${\mathcal P}$ lying over $b$.

\begin{definition} The set ${\mathcal P}$ is said to be {\em $n$-correlated}
if there is a proper Zariski-closed subset $F_n \subset X_B^n$ such that
${\mathcal P}_B^n \subset F_n$. 
\end{definition}

For instance, a subset ${\mathcal P}$ is 1-correlated if and only if it is not
Zariski-dense in $X$; in which case it is easy to see that over some nonempty
open subset $U$ 
in $B$ the number of points of ${\mathcal P}$ in each fiber is uniformly
bounded.  This simple observation is generalized by the following lemma.

\begin{lemma}[cf. Lemma 1.1 in \cite{C-H-M}, or Lemma 1 in \cite{A:stintegral}]
\label{Lem:correlation}
Let $X
\rightarrow B$ a projective family of smooth irreducible curves, ${\mathcal P}
\subset X(K)$ an
$n$-correlated subset.  Then there exists a nonempty open subset $U \subset B$
and 
an integer $N$ such that for every $b \in U$ we have $\# {\mathcal P}_b \leq
N$. 
\end{lemma} 

We find it instructive to include a short  proof of this simple lemma, in which
it will be clear 
that the argument only works for a family of curves. A modification which does
work in higher dimension will be discussed later in this article - see
Sections \ref{Ssec:unif-higher} and \ref{Ssec:prolongable}. 

 Let $F_n = \overline{{\mathcal P}_B^n}$ be the Zariski
closure, $U_n = X_B^n \setmin F_n$ the complement.  We now define Zariski-open
and Zariski-closed subsets $U_{i-1}$ and $ F_{i-1} \subset X_B^{i-1}$ by
descending 
induction as follows: we take
 $U_{i-1} = \tau_i(U_i)$, and set  $F_{i-1} = X_B^{i-1} \setmin U_{i-1}$ to be
the complement.  

Note that  over $U_{i-1}$ the map $\tau_i$ restricts to a finite map on $F_i$. 
In fact, by definition, if $x \in U_{i-1}$ then $\tau_i^{-1}(x) \not\subset
F_i$ and hence $\tau_i^{-1}(x) \cap F_i$ is a finite set, since
$\tau_i^{-1}(x)$ is a curve.  Thus there exists
$d_i \in {\mathbb N}$ 
such that
$$\# \tau_i^{-1}(x) \cap F_i \leq d_i \text{\ for\ }x \in U_{i-1}.$$
Let $U = U_0 \subset B$.  We claim that over $U$ the number of points of
${\mathcal P}$ in each fiber is bounded.  Consider a point $b \in U$.

Case 1: ${\mathcal P}_b \subset F_1$.  

In this case, we have $\# {\mathcal P}_b \leq d_1$.

Case 2: There exists some $Q \in {\mathcal P}_b \setmin F_1$ where
$X_Q
\cap {\mathcal P}_b^2 \subset F_2$.

In this case, we have $\# {\mathcal P}_b \leq d_2$.

Case $i$: There exists $Q \in (P_1, \cdot\cdot\cdot, P_{i-1}) \in {\mathcal
P}_b^{i-1} \setmin F_{i-1}$ where $X_Q \cap {\mathcal P}_b^i \subset F_i$.

In this case, we have $\# {\mathcal P}_b \leq d_i$.

As $X_Q \cap {\mathcal P}_b^n \subset F_n$ for all $Q \in {\mathcal P}_b^{n-1}$ by definition, the cases will
be exhausted at some stage when $i \leq n$.  Thus
$$\# {\mathcal P}_b \leq \max \{d_i\}.$$  \qed

\subsection{Lang's conjecture and correlation}
The remarkable observation of \cite{C-H-M} is that Lang's conjecture 
implies that the set of $K$-rational points on curves of genus $g>1$ is
$n$-correlated, for sufficiently large $n \in {\mathbb N}$.  

\begin{proposition}[cf. Lemma 1.1 in  \cite{C-H-M}]\label{Prop:curves} Let $X
\rightarrow B$ be a 
projective family of smooth irreducible curves of genus $g \geq 2$ over a
number field
$K$.  Assume that Lang's conjecture holds true.  Then
$X(K)$ is $n$-correlated for sufficiently large $n \in {\mathbb N}$.
\end{proposition}

In order to deduce the uniformity assertion of  Section \ref{Sec:uniformity}
from Proposition \ref{Prop:curves}, 
one starts with $X \rightarrow B$, a ``comprehensive''
projective family of smooth irreducible curves of genus $g \geq 2$ in which
``all curves appear'', namely,
for any projective smooth curve $C$ of genus $g$ defined over $K$ there exists
a morphism $Spec K \rightarrow B$ satisfying
$$\begin{CD}
C \simeq \text{Spec} K \times X @>>> X \\
@VVV @VVV \\
\text{Spec} K @>>> B\ . \\
\end{CD}$$ 
Such a family always exists over a suitable Hilbert scheme, since all curves of
genus $g>1$ are canonically polarized.
We set ${\mathcal P} = X(K)$.  There exists a nonempty Zariski-open subset $U_0
\subset 
B$ and an integer $N_0$ such that $\# {\mathcal P}_b \leq N_0$ for $b \in U_0$
by 
Lemma \ref{Lem:correlation}.  Now take $B_1 = B \setmin U_0$, and  apply the lemma to
the family $X_1 
= X  \times_B B_1
\rightarrow B_1$ to obtain a new nonempty Zariski-open subset $U_1 \subset
B_1$ and an integer $N_1$ such that $\# {\mathcal P}_b \leq N_1$ for $b \in
U_1$, and so on.  By noetherian induction, we have the uniformity
assertion.

\subsection{The fibered power theorem}
It is easy to see that Proposition \ref{Prop:curves} follows from the 
{\em Fibered Power Theorem:}

\begin{theorem}\label{Th:fibered-power} Let $\pi:X \rightarrow B$ be  
a projective family of varieties of general type, with $B$ irreducible, defined
over a field 
$K$ of characteristic 0.  Then there exists a positive integer $n$,
 a variety of general type $W_n$ over $K$ with $\dim W_n > 0$, and 
 a dominant rational map
$$r:X_B^n \dashrightarrow W_n.$$
\end{theorem}

To see how Proposition \ref{Prop:curves} follows from this theorem, first note
that we may replace $B$ by an irreducible 
component. Next, by Lang's conjecture, there exists a
proper Zariski-closed subset  $G_n$ of $W_n$ which contains all the
$K$-rational  
points of $W_n$. Let $(X_B^n)_{dom}$ be the domain of the rational
map $X_B^n \dashrightarrow W_n.$  Denote   ${\mathcal P} = X(K)$.  Then for 
a point $P \in {\mathcal P}_B^n \subset X_B^n$, we
have 
either $P \in X_B^n \setmin (X_B^n)_{dom}$ or $P \in r^{-1}(G_n)$.   Then we
only 
have to set  $$F_n =  [X_B^n \setmin (X_B^n)_{dom}] \cup  [r^{-1}(G_n)].$$

In the case of curves, the Fibered Power Theorem was proven in \cite{C-H-M}
Theorem 1.3, 
using the following observation: let $X \to B$ be a family of curves of genus
$g>1$ as above.  For each $n$ we have a rational map $X^n_B \das
\bM_{g,n}$. Denote by $W_n$ the image of this map. Then an argument is given in
\cite{C-H-M} which in effect proves the following:

\begin{quote} \em 
For large $n$, the image variety $W_n$ is a variety of general type. 
\end{quote}

The argument in  \cite{C-H-M} uses the compactification $\obM_{g}$, the
moduli space of {\em stable} curves, in an essential way. A similar argument
was used 
in \cite{Hassett:surfaces},  Theorem 1,  for the case of surfaces. This is
precisely the 
line of proof 
we will take in this paper for abelian varieties (see Theorem \ref{Th:lgt}). In
higher dimension, it was necessary to give a different argument in
\cite{A:fibered}, Theorem 0.1, since complete moduli spaces of stable  
varieties in dimension $>2$ are not known to exist in general. All proofs use
deep results about {\em weak positivity} of the push forward of pluricanonical
sheaves; in the present paper we use such a result due to Kawamata
\cite{Kawamata:MMandKod} about
abelian varieties.

\subsection{Uniformity in higher dimension}\label{Ssec:unif-higher}

% As we have noted,  the arguments given above need to be modified
%in order to 
%apply  to varieties of dimension $\geq 2$. We now provide such a modification,
% following  \cite{A:quad}, \cite{Pacelli:rational}, and
% \cite{A-Voloch:fibered}. 

%

In order to prove a uniformity result in dimension $>1$, one needs to
modify the {\em statement} appropriately, and then adjust the
proof.

First, as a variety of general type may contain a subvariety which is 
not of general type, on which there may be infinitely many $K$-rational points,
we have to modify the uniformity statement (this issue does not come up in this
paper). Such a subvariety is called {\em 
exceptional}. Given a family $X \to B$ of varieties of general type, it is
natural to restrict attention to points in $X(K)$ which do 
not lie on any exceptional subvariety, which we denote $X(K)^{nex}$ (``nex''
for non-exceptional). We arrive at the following  statement (see
\cite{A-Voloch:fibered}, Theorem 1.5): 

\begin{quote}  Assume Lang's conjecture holds true. Fix a family of varieties
of general type $X \to B$. Then there 
exists a number
$N(X \to B,K)$ depending only on the given family and the number field
$K$ such that  for any $b\in B(K)$, the number of non-exceptional points on
the fiber $X_b$ is bounded:
$$\# X(K)^{nex} < N(X \to B,K).$$
\end{quote}

Second, we still need to modify Lemma \ref{Lem:correlation} for higher
dimension.
Such a modification was given in \cite{A:quad}, \cite{Pacelli:rational}, and
\cite{A-Voloch:fibered}, and is essential in this paper as well. We will
discuss some aspects of it in the course of the proof of the Main Theorem (see
Section \ref{Ssec:prolongable}).  First, we would like to adjust the notion of
correlation for the higher dimensional case:

\begin{definition}\label{Def:strong-correlation}
 Let $\tau:X \to B$ be a projective surjective morphism of
reduced schemes of finite type over
a  field $K$,  with $B$ irreducible.
 Denote the dimension of the
generic fiber of $\tau:X \to B$ by $d$. Fix a subset
$\cQ\subset X(K)$, and denote by $G_k$ the Zariski closure of $\cQ_B^k$ in
the fibered power $X_B^k$.  
 We say that $\cQ$ is {\em strongly
$k$-correlated} with respect to  $\tau:X \to B$, if 
every irreducible component of $G_k$ which dominates $B$ has relative dimension
$<kd$ over $B$. 
\end{definition}

We will now see how to reduce a question of strong correlation
to a question of correlation.
Let $X = X_1 \cup \ldots \cup X_c$ be a decomposition into irreducible
components. Let $X_i'\to X_i$ be the normalization, and let $X'_i \to B'_i \to
B$ be the Stein factorization. There is a dense open set $U_i\subset X_i$ over
which  $X_i'\to X$ is an isomorphism. Therefore the set $\cQ_i' = U_i \cap
\cQ$ sits 
naturally in $X_i'(K)$.

\begin{proposition}\label{Prop:strong-correlation} Assume that for each $X_i$
of relative dimension $d$ over 
$B$, the set $\cQ_i'$ is $k_i'$-correlated with respect to $X_i'\to B_i'$. Then
for large $k$, the set $\cQ$ is strongly $k$-correlated for $X \to B$. 
\end{proposition}

{\bf Proof.} We make a number of reduction steps leading to the proposition.

\subsubsection{}\label{Sssec:strong-cor-open}
{\em  Let $B'\subset B$ be a nonempty open set, $X'\subset \tau^{-1}
B'$ a dense open set, and 
$\cQ' = \cQ \cap X'$. Assume $\cQ'$ is strongly $k$-correlated. Then $\cQ$ is 
strongly $k$-correlated as well.}

This is immediate from the definition.

In particular, we may as well assume that $B$ is normal and  $\tau:X \to B$ is
flat. We may also replace $X$ by a birational modification,
since we can restrict to the points of $\cQ$ where this modification is an
isomoprhism. 

\subsubsection{} {\em Let $X_1,\ldots,X_{c'}$ be the irreducible components of
$X$ of relative dimension $d$ over $B$. Let $X' = X_1 \cup\ldots\cup X_{c'}$
and 
let $\cQ' = \cQ \cap X'$. If $\cQ'$ is strongly $k$-correlated then so is
$\cQ$.  } 

One can write $$X^k_B = \mathop{\bigcup}\limits_{1 \leq i_j \leq c} X_{i_1}
\times_B \ldots \times_B X_{i_k}.$$ 

If some component $X_i$ has relative dimension $<d$, then any term
involving $X_i$ in $X^k_B = \cup_{1 \leq i_j \leq c} X_{i_1} \times_B
\ldots \times_B X_{i_k}$ has relative dimension $<kd$.

Thus we may as well assume that all fibers of $X \to B$ have pure dimension
$d$. 

\subsubsection{} {\em Let $\cQ_i=X_i\cap \cQ$. Assume $\cQ_i$ is strongly
$l$-correlated. Then for $k=l\cdot c$, $\cQ$ is strongly $k$-correlated.}

 By the box principle, every term in
the expression $$X^k_B = \mathop{\bigcup}\limits_{1 \leq i_j \leq c} X_{i_1}
\times_B  
\ldots \times_B X_{i_k}$$ has at least one $i_j$ appearing at least $l$ times. 
Considering the projection on those factors, it follows that $\cQ$ is strongly
$k$-correlated. 

\subsubsection{}  {\em Denote by $c_i$ the degree of $B_i'\to B$. Assume
$\cQ_i'$  in the proposition is $k_i'$-correlated with respect to $X_i' \to
B_i'$ .   Then for $k = c_ik_i'$
we have that $\cQ_i$ is strongly $k$-correlated with respect to
$X_i \to B$.}   

Let $G$ be an irreducible component of ${X'_i}^k_B$. A point on $G$ corresponds
to a $k$-tuple of points on a fiber of $X'_i$ over $B$, which fall into the
$c_i$ 
different components of this fiber, which are identified as fibers of $X'_i$
over 
$B'_i$. By the box principle there is a subset $J\subset \{1,\ldots k\}$ of
size at least $k/c_i = k_i'$ such that for a point $(P_1,\ldots P_k)\in G$, all
the 
$P_j$ for $j\in J$ lie in the same component of the fiber. In other words,  the
projection  $G\to (X_i)^{J}_B$ to the factors in $J$ maps
$G$ onto the closed subset $(X_i)^{J}_{B'_i}$.   Since $\cQ_i'$  is
$k_i'$-correlated for $X_i'\to B_i'$, we have that $(\cQ_i')^k_B\cap G$ 
is not dense in $G$, which implies that  $\cQ_i'$ is strongly $k$-correlated
with respect to $X_i' \to B$. Step \ref{Sssec:strong-cor-open} in the proof
implies that $\cQ_i$ 
is strongly $k$-correlated as well. \qed

\section{Moduli of stable quasi-abelian pairs}\label{Sec:Alex}

Both statement and proof of our main theorem depend on the existence
of a good compactification of the moduli space of abelian varieties. We now
review some essential facts about these spaces which we will utilize.

\subsection{Abelian schemes} Recall that a principally polarized abelian scheme
$(A\to S, \lambda)$ is an abelian scheme $A \to S$ (with a zero section $S\to
A$) and an {\em isomorphism} $\lambda:A \to \Pic^0(A/S)$ which {\em locally
over $S$} is a polarization induced by a relatively ample invertible 
sheaf.  

The moduli category $\fA_g$ of  principally polarized abelian schemes
(morphisms given by fiber diagrams) is a
Deligne-Mumford stack. It admits a coarse moduli scheme $\bA_g$, which is 
quasi-projective over $\Spec \ZZ$ (see
\cite{Mumford:GIT}, \cite{Faltings-Chai}). In analogy with the moduli spaces of
curves, one would like to have a good compactification of $\fA_g$ in a
canonical way, and possibly also an analogue of $\ocM_{g,n}$. 

Beginning with \cite{A-M-R-T} and through the work of many authors (see
\cite{Faltings-Chai}), an infinite collection of  troroidal compactifications
``$\ofA_g$'' of 
$\fA_g$ were constructed, depending on choices of ``cone decompositions''. In
general these compactifications are not moduli stacks of any explicitly
described families of ``stable objects''. It is however shown in
\cite{Faltings-Chai} that each of these compactifications carries a family of
semiabelian varieies. If one takes the formal completion of such a $\ofA_g$ at
a point, one can apply {\em Mumford's construction} (See 
\cite{Mumford:degeneration}, \cite{Faltings-Chai}) and get a toroidal 
compactification of the family of semiabelian schemes, but this
compactification again depends on ``degeneration data'', and one has a
serious problem in gluing these together.

 These issues were recently resolved in the work of Alexeev and Nakamura
\cite{Alexeev-Nakamura}, \cite{Alexeev:SQAV}, \cite{Alexeev:toric}. See also
the related \cite{Nakamura:level}.

\subsection{Stable quasi-abelian pairs}
A first important step is to
change the original moduli problem in a way which surprizingly simplifies the
situation. Instead of working with  principally polarized abelian schemes, one
forgets the zero section of the abelian schemes and instead one insists that
the polarization come 
from a global relatively ample invertible sheaf; in fact, since we work with
principal polarizations, this sheaf has a unique divisor $\Theta$. To this end,
a smooth principally polarized {\em  quasi-abelian scheme} $(P \to S,\Theta)$
is 
a {\em 
torsor} $P \to S$ on an abelian scheme $A \to S$, and a relatively ample
divisor $\Theta\subset P$ which behaves like a principal polarization, in the
sense that its Hilbert polynomial is $n^g$ (see \cite{Alexeev:SQAV}, 1.4). The
moduli 
stack of smooth principally polarized 
quasi-abelian schemes is canonically isomorphic to $\fA_g$
(\cite{Alexeev:SQAV}, Theorem 1.15). It admits a 
universal family which we denote $(\fA_{g,1}\to
\fA_g,\tilde{\boldsymbol{\varTheta}})$;
the added subscript $1$ indicates that $\fA_{g,1}$ is the moduli stack for
smooth principally 
polarized  quasi-abelian schemes {\em with one marked point}. (It should
not be 
confused with Mumford's notation for level structure).

Next, Alexeev and Nakamura make a 
canonical choice of degeneration data in Mumford's construction. Over a
discrete valuation ring, this gives a
canonical way to compactify a torsor on a semiabelian scheme with smooth
principally 
polarized  quasi-abelian generic fiber. The  fibers $(P,\Theta_P)$ of
this construction are called {\em stable principally polarized 
quasi-abelian varieties}, and they can be characterized explicitly, see
\cite{Alexeev:SQAV}, Definition 1.11. (The reader is advised not to be confused
by this nomenclature: a smooth quasi-abelian variety is by definition also
a stable quasi-abelian variety.)

Finally, it follows from Alexeev's work \cite{Alexeev:SQAV},
\cite{Alexeev:toric} that the category 
of  stable principally polarized  quasi-abelian schemes is a Deligne-Mumford
stack $\ofA_g$ admitting   a projective coarse moduli scheme
$\obA_g\to \Spec \ZZ$.   On
the level of geometric points, Alexeev shows that $\ofA_g$ agrees with the
so called ``Second Voronoi Compactification'' $\ofA_g'$, which is a very
special toroidal compactification of $\fA_g$. Indeed, there is  a 
morphism $\ofA_g'\to\ofA_g$ which is one-to-one on geometric points
(\cite{Alexeev:SQAV}, Theorem 4.5). Alexeev remarks that in general this is
not an isomorphism.

Again, we
denote the universal family by $(\ofA_{g,1},\ofTheta) \to \ofA_g$. 
We denote by $\ofA_{g,n}$ the fibered power $(\ofA_{g,1})^n_{\ofA_g}$. (This
is, to some extent, in analogy with the space of stable pointed curves
$\ocM_{g,n}$, although we do not use Knudsen's stabilization.) Denote
by $p_i:\ofA_{g,n}\to\ofA_{g,1}$ the projection to the $i$-th factor. We have
a natural relatively ample divisor $\ofTheta_n \subset \ofA_{g,n}$ defined by
$\ofTheta_n=\sum_i p_i^*\ofTheta$. 

We denote by $\obA_{g,n}$ the coarse moduli spaces of $\ofA_{g,n}$, and
by $\obTheta_n \subset \obA_{g,n}$ the  image of $\ofTheta_n$. A-priori these
are Artin algebraic spaces (see \cite{Keel-Mori}), but since some multiple
$m\ofTheta_n$ descends to a Cartier divisor on $\obA_{g,n}$ and is relatively
ample, these are projective schemes over $\Spec \ZZ$. 

\subsection{Properties of stable pairs} 

We now collect a few properties of stable principally polarized  quasi-abelian
schemes, which we will use in the next section.

To save words, we will refer to a stable principally polarized
quasi-abelian 
scheme $(P, \Theta)$ (always assumed flat over a base scheme $S$) as a {\em
stable pair}. 

The first two items are included in \cite{Alexeev:SQAV}, Definition 1.11 and
\cite{Alexeev:toric}.

\subsubsection{} For a stable pair $(P,\Theta)$ over a field, the underlying
stable quasi-abelian variety $P$ is proper and reduced,  and
$\Theta$ is an ample Cartier divisor.

%\subsubsection{}  The pair $(P,\Theta)$ has
%semi-log-canonical singularities. 
%
%See \cite{Alexeev:pairs}, Theorem 3.10.  

\subsubsection{} Let  $(P\to S,\Theta)$ be a stable pair over $S$. Let
$P_0\subset P$ be an open subset, consisting of exactly one irreducible
component of the smooth locus in every fiber. Then $A =
\Aut^0(P_0/S)$ is semiabelian, $P_0$ is an $A$-torsor, and $A$ independent of
the choice of $P_0$. Over a 
field, $P$ is stratified by finitely many orbits of $A$.

\subsubsection{}\label{Sssec:sqav-toroidal} Let $U_S \subset S$ be a toroidal
embedding over a field, and 
let $(\pi:P\to S,\Theta)$ be a stable pair over $S$, such that $P \to S$ is
smooth over the open set $U_S$. Let $U_P = \pi^{-1} U_S$. then $U_P \subset P$
is a toroidal  embedding, and $P \to S$ is a toroidal morphism. (Toroidal
morphisms are defined in \cite{A-Karu}, Definition 1.2.)

Indeed, Mumford's construction, as well as Alexeev's variant in
\cite{Alexeev:toric}, is by definition toroidal!

(If one is working in mixed characteristics, one only needs to replace
``toroidal'' by ``log-smooth'' in the sense of K. Kato.)

\subsubsection{}\label{Sssec:sqav-powers} Let $(P, \Theta)$ be a stable pair
over a field. Let $p_i:P^n 
\to P$ be the projection onto the $i$-th factor, and consider the divisor
$\Theta_n = \sum_ip_i^*\Theta$.  Then $(P^n, \Theta_n)$ is a stable pair.

This follows immediately from \cite{Alexeev:SQAV}, Definitions 1.10, 1.11, or
\cite{Alexeev:toric}. 

\subsubsection{}\label{Sssec:sqav-gorenstein} In the situation of
\ref{Sssec:sqav-toroidal},  suppose $S$ 
has Gorenstein singularities. Then the scheme $P$ has Gorenstein
singularities.  

This is a general fact about toroidal morphisms with {\em reduced fibers} and
{\em no horizontal divisors.} The proof is easy using the associated polyhedral
complexes,  see \cite{A-Karu}, Lemma 6.1.  

We note that, since a toroidal embedding has rational singularities, it follows
that $P$ has rational Gorenstein singularities, hence {\em canonical}
singularities.  
This is a refinement of Alexeev's \cite{Alexeev:pairs}, Lemma 3.8.

\subsubsection{}\label{Sssec:sqav-lc} Suppose the base field has
characteristic 0.  In the situation of
\ref{Sssec:sqav-gorenstein}, the pair $(P, 
\Theta)$ has log-canonical singularities.

Indeed, the proof of \cite{Alexeev:pairs}, Theorem 3.10 applies word-for-word,
if we do not add the central fiber  $P_0$ (this only makes the proof simpler).

%\smallskip

Finally, we have the following crucial extension property, proved in
\cite{Alexeev:SQAV} and \cite{Alexeev:toric}:
\subsubsection{} Let $S$ be the spectrum of a discrete valuation ring $R$, with
generic point $\eta$. Let $(P_\eta ,
\Theta_\eta)$ be a stable pair. Then there exists a finite separable extension
of discrete valuation rings $R \subset R_1$, with spectrum $S_1$ and generic
point $\eta_1$, and a stable pair $(P_1 \to S_1, \Theta_1)$ extending
$(P_{\eta_1} , \Theta_{\eta_1})$.

This result immediately extends to dedekind domains.

We call $(P_1 \to S_1, \Theta_1)$ a {\em  stable quasi-abelian model} of
$(P_{\eta_1}, \Theta_{\eta_1})$.

\subsection{Relation with the N\'eron model}
Let $R$ be a discrete valuation ring, $S =\Spec R$, with generic point $\eta$
and special point $\fp$. Let $(P\to S, \Theta)$ be a stable pair, and assume
$A = P_\eta$ is smooth. For simplicity we also assume there is a section
$s:S 
\to 
P$ landing in the smooth locus of $P\to S$. This makes $A$ into an abelian
variety. Denote by $\cN_A \to S$ the N\'eron model of $A$.

\begin{proposition}\label{Prop:sqav-neron} There is a unique morphism
$$f:{\mathcal N}_A \rightarrow {P}$$
extending the isomorphism
$$({\mathcal N}_A)_{\eta}
\overset{\sim}{\longrightarrow} A.$$
\end{proposition}

{\bf Proof.} The formation of  $\cN_A \to S$ commutes with \'etale base change.
 Once we  prove the proposition after such a base change, the uniqueness
 implies that we can descend back to $S$. Thus we may replace $R$ by its strict
 henselization. 

Now by construction, the  stable quasi-abelian model $({P}\to
S,\tilde{\varTheta}) $ 
can be viewed as a ``compactification'' of a semi-abelian scheme ${\mathcal
A}_{0} = \Aut^0(P/S)
\rightarrow S$ with the origin identified with the section $s$.  Note that the
construction of 
the  stable  quasi-abelian model gives an action 
$${\mathcal A}_{0} \times {P} \rightarrow {P}$$
extending the addition law on ${\mathcal A}_{0}$.

Note that there is also a natural inclusion ${\mathcal A}_{0}
\hookrightarrow 
{\mathcal N}_A$ as the zero component.

Denote by
$$M_i, \quad i = 1, \cdots, t$$
the components of the fiber $({\mathcal N}_A)_\fp$ of the N\'eron model over
$\fp$.  For each 
$i$, we have an open neighborhood
$${\mathcal N}_i = {\mathcal N}_A \setmin \cup_{j \neq i}M_j.$$
We may choose the numbering so that ${\mathcal N}_1 = {\mathcal A}_{0}.$
We have
$${\mathcal N}_A = \cup_i {\mathcal N}_i$$
and 
$${\mathcal N}_i \cap {\mathcal N}_j = ({\mathcal N}_A)_{\eta} \quad \forall
i\neq j.$$ 

Since $R$ is strictly henselian, we can choose, for each $i$,
a section $$s_i:S \rightarrow {\mathcal N}_i$$
such that $s_1 = s$.
 The schemes ${\mathcal N}_i$ can be viewed as ${\mathcal 
A}_{0}$-torsors, and the choice of the $s_i$ gives a trivialization of these
torsors. 

We denote by
$$(s_i)_{\eta}, \quad i = 1, \cdots, t$$
the corresponding rational points on $({\mathcal N}_A)_{\eta} = {P}_{\eta}$.
Since  
${P} \rightarrow S$ is proper, we can extend $(s_i)_{\eta}$ to 
sections
$$t_i\in P(S)$$
such that
$$(t_i)_{\eta} = (s_i)_{\eta}.$$

We define
$$f_i:{\mathcal N}_i \rightarrow {P}, \quad i = 1,
\cdots, t$$ 
as follows.  Given a scheme $T$ over $S$ and a
point $z \in 
{\mathcal N}_i(T)$ there exists a unique point $a\in {\mathcal A}_{0}(T)$ such
that $z = a \cdot s_i$.  Here the notation
$a \cdot s_i$ stands for the action ${\mathcal A}_{0} \times {\mathcal N}_i
\rightarrow {\mathcal 
N}_i$.  Define $f_i(z) = a \cdot t_i$.  Here the notation $a \cdot t_i$ stands
for the 
action ${\mathcal A}_{0} \times P \rightarrow P$. This
rule is clearly 
functorial and therefore defines a morphism.

We claim that the morphisms $f_i$ are independent of the choice of $s_i$ and
coincide 
on $({\mathcal N}_A)_{\eta}$.  In fact, given $s_i' \in {\mathcal N}_i(R)$
consider the corresponding
$t_i'$ and $f_i'$.  There exists $b_i \in {\mathcal A}_{0}$ such
that $s_i = b_i 
\cdot s_i'$.  Therefore, $(t_i)_{\eta} = (b_i)_{\eta} \cdot (t_i')_{\eta}$,
which implies 
that $t_i = b_i \cdot t_i'$, since $P$ is separated.  Therefore, we conclude  
\begin{eqnarray*}
f_i(z) &=& f_i(a \cdot s_i) = a \cdot t_i = a \cdot (b_i \cdot t_i') = (ab_i)
\cdot t_i' \\ 
&=& f_i'((ab_i) \cdot s_i') = f_i'(a \cdot (b_i \cdot s_i')) = f_i'(a \cdot
s_i) = f_i'(z).
\end{eqnarray*} 
The same argument shows that given a scheme $T$ over
$\eta$ and a point $z \in {\mathcal N}_A(T)$ we have $f_i(z) = f_j(z)$, i.e.,
the 
morphisms $f_i$ coincide over the intersection of their domain $({\mathcal
N}_A)_{\eta}$. 

Since ${\mathcal N}_A$ is covered by the ${\mathcal N}_i$ and ${\mathcal N}_i
\cap {\mathcal N}_j = 
({\mathcal N}_A)_{\eta}$ whenever $i \neq j$, it follows that the $f_i$ glue
together to 
give a morphism ${\mathcal N}_A \rightarrow {\mathcal A}$ as required in the
Proposition. \qed

\subsection{Tautological families over moduli spaces} Since $\obA_g$ is not a
fine 
moduli scheme, it does not have a universal family. Our arguments below
depend on the geometry of families, therefore it is useful to have some
approximation of a universal family, which, following \cite{C-H-M}, one calles
{\em a tautological family}. We need such a family with a strong equivariance
property. This is summarized in the following statement. 

\begin{proposition}\label{Prop:taut}
Let $W_0 \subset \obA_g$ be a closed integral subscheme, and let  $W_1 \subset
\obA_{g,1}$ be the reduced scheme underlying its inverse image. Then there
exists a 
 projective, normal  integral scheme $B$, and a family of $g$-dimensional
stable 
pairs $(A \to B, \Theta)$, satisfying the following properties:

\begin{enumerate}
\item The natural moduli morphisms $B \to  \obA_g$ and $A \to \obA_{g,1}$
are finite and generically \'etale, with images $W_0$ and $W_1$, respectively. 
\item Denote $G = \Aut_{\obA_g}(A \to B, \Theta)$, namely the group of
automorphisms of  the gadget $(A \to B, \Theta)$ commuting with the morphism $B
\to \obA_g$. Then the two morphisms $B / G \to \obA_g$ and
$A/G\to \obA_{g,1}$ are birational onto their images.
\end{enumerate}
Moreover, in this situation consider the moduli morphism $A^n_B \to
\obA_{g,n}$. There is a diagonal action of $G$ on $A^n_B$, and $A^n_B/G\to
\obA_{g,n}$ is again birational onto the image.
\end{proposition}

{\bf Proof.} The existence of a family $(A \to B, \Theta)$ satisfying  
condition 1 is an immediate consequence of  \cite{Kollar:projectivity},
Proposition 
2.7. (Koll\'ar attributes  the proof to M. Artin. See also discussion in
\cite{C-H-M}. Proofs of this fact have been given  by a number of authors
through the years.)

We now wish to replace this family by one which satisfies the equivariance
condition 2. First, we may assume that the function field extension $K(W_0)
\subset K(B)$ is Galois, by going to the Galois closure. Denote the Galois
group $G_0$. Write $\eta$ for the generic point of $B$. Second, we may assume
that the geometric points of the finite group  $H=\Aut_\eta(A_\eta,
\Theta_\eta)$ 
are all rational over $K(B)$ - simply pass to a suitable finite extension and
take Galois closure again. Since the $\Isom$ scheme is proper over $B$, the
automorphisms of  $(A_\eta, \Theta_\eta)$ extend to automorphisms  of the
family $(A\to B,\Theta)$ over $B$. Third, in a similar manner we can ensure
that for any $g\in G_0$ we have a $B$-isomorphism $(A\to
B,\Theta)\stackrel{\sim}{\longrightarrow} g^*(A\to B,\Theta)$. Now consider
the set $$G = \{(g,h) | g\in G_0, h:(A\to
B,\Theta)\stackrel{\sim}{\longrightarrow} g^*(A\to B,\Theta)\}.$$

Note that for any $(g, h)\in G$ and $g'\in G_0$  we can define
$$h^{\{g'\}}:{g'}^*(A\to  B,\Theta)\stackrel{\sim}{\longrightarrow}
{g'}^*g^*(A\to B,\Theta)\}.$$ We can 
now define 
$$(g,h)(g',h') = (gg',h^{\{g'\}}h').$$ 
We leave it to the reader to check that this is a group. It is now immediate to
verify condition (2) and the fact that $A^n_B/G\to
\obA_{g,n}$ is  birational onto the image.

\section{Proof of the Main Theorem}\label{Sec:proof}

\subsection{Stably $S$-integral points}
Let $A$ be a principally polarized abelian variety, with  theta divisor
$\Theta$, defined over a number field $K$. Let $S$ be  a finite number of
places in $K$, and denote by ${\mathcal O}_{K,S}$ the ring of $S$-integers.
For an extension $L$ of $K$, we denote by $S_L$ the set of places in $L$ over
$S$. 

\begin{definition}\label{Def:stably-integral}
A $K$-rational point $P \in (A \setmin \Theta)(K)$ is called a {\em
stably $S$-integral point} if {\em there exists} a finite extension $L \supset
K$, 
and a
 stable quasi-abelian model $(\cP\to \Spec \cO_{L,S_L},\cTheta)$,
such that 
$P \in (\cP \setmin \cTheta)({\cO}_{L,S_L})$, namely $P$ is integral
on 
the complement of $\cTheta$ in $\cP$.
\end{definition}

The following result shows that the existential quantifier is not important in
the definition: 

\begin{proposition}\label{Prop:stint-any-L}
Suppose $P \in (A \setmin \Theta)(K)$ is
stably $S$-integral. 

Let $L' \supset K$ be {\em any} finite extension for which there exsits a 
 stable  quasi-abelian model $(\cP'\to \Spec \cO_{L',S_{L'}},\cTheta')$. Then
$P \in (\cP \setmin \cTheta')({\mathcal O}_{{L'},S_{L'}})$.
\end{proposition}

{\bf Proof.}   Let $L \supset K$ be a finite extension satisfying the
conditions in the definition.
Take a Galois extension $M \supset L'$ which contains $L$, with Galois
group $G = Gal(M/L')$.  Then by the functoriality of the stable model
\begin{eqnarray*}
P &\in &\{(\cP \setmin \cTheta)({\mathcal
O}_{L,S_{L}}) \mathop{\otimes}\limits_{L} M\} \cap \{(A \setmin \Theta)(K)
\mathop{\otimes}\limits_{K} M\} \subset 
(\cP_M \setmin \cTheta_M)({\mathcal O}_{M,S_M})^G \\
&=& (\cP \setmin \cTheta)({\mathcal O}_{L',S_{L'}}).
\end{eqnarray*}
This completes the proof.

Stably integral points have a nice characterization in terms of moduli:

\begin{proposition}
Let $P \in (A \setmin \Theta)(K)$. Consider the associated moduli morphism
$P_m:\Spec K \to A \to \obA_{g,1}$. Then $P$ is  
stably $S$-integral if and only if $P_m$ is an $S$-integral point on $
\obA_{g,1} \setmin \obTheta$. 
\end{proposition}

{\bf Proof.} Clearly $P_m$ is a rational point on $\obA_{g,1}
\setmin \obTheta$, so to check that it is $S$ integral we may pass to a finite
extension field. Let $L\supset K$ be an extension such that there exists a
stable quasi-abelian model $(\cP\to
\Spec\cO_{L,S_L},\tilde{\varTheta})$. Since $\cP$ and $\obA_{g,1}$ are proper,
we have  morphisms 
$\tilde{P}:\Spec\cO_{L,S_L} \to \cP$ and $\tilde{P_m}:\Spec\cO_{L,S_L} \to
\obA_{g,1}$. Note that by the coarse moduli property, $\tilde{\varTheta}$ is
the set theoretic inverse image of $\obTheta$. Then $P$ is stably  $S$-integral
if and only if $\tilde{P}$ is disjoint  from $\tilde{\varTheta}$, if and only
if $\tilde{P_m}$ is disjoint from $\obTheta$. \qed

\subsection{Reduction to moduli}\label{Ssec:red-to-mod} Fix a 
 family $(\fA\to \fB, \fXi)$ of smooth principally polarized quasi-abelian
varieties with 
theta divisors. We say that a point
${P} \in ({\fA} \setmin 
\fXi)(K)$ is  {\em 
stably  
$S$-integral}, if 
it is stably $S$-integral on the fiber of $\fA \to \fB$ on which it lies.

Let ${\mathcal P} \subset ({\fA} \setmin \fXi)(K)$  be the set
of stably 
$S$-integral points. 
\begin{proposition}\label{Prop:correlation} 
Assume that the Lang-Vojta conjecture holds true. Then for a sufficiently large
integer $n$, the set  ${\mathcal P}$ is $n$-correlated.
\end{proposition}
 
{\bf Proof.}
We may assume
$\fB$ is irreducible (by taking irreducuble components one by one) and 
hence is   a variety.

Following \cite{C-H-M}, we would like to reduce the situation to a situation on
moduli spaces.
There is the natural moduli morphism $\nu:\fB \rightarrow \obA_g$ from $\fB$ to
the 
coarse moduli space of stable pairs.  We denote
by $W_{0}$ the image $\nu(\fB)$ under this map. There is also a compatible
dominant morphism $\fA  \to W_1 \subset \obA_{g,1}$, creating a commutative
diagram: 

$$\begin{array}{ccccc}  \fA & \to&W_1 & \subset&  \obA_{g,1}\\
\dar & &\dar & & \dar\\
			\fB &\rightarrow&W_0 & \subset&  \obA_g.
\end{array}$$

Recall that we have characterized stably integral points in terms of their
image in moduli. Thus the Proposition follows immediately from the following
purely geometric result:

\begin{theorem}\label{Th:lgt}
Let $W_0\subset (\obA_g)_{\CC}$ be  closed subvariety, and suppose $W_0\cap
(\bA_g)_{\CC} \neq \emptyset$ (thus the generic point of $W_0$ parametrizes a
smooth quasi-abelian variety).  Let $W_n \subset (\obA_{g,n})_{\CC}$ be the
reduced scheme underlying  the inverse image, and $\Theta_{W_n} = W_n \cap
\obTheta_n$. Then 
for large integer 
$n$, the pair $(W_n, \Theta_{W_n})$ is of logarithmic general type.  
\end{theorem}

In view of Proposition \ref{Prop:taut}, it suffices to prove the following:

\begin{proposition}\label{Prop:quotients} Let $({ P}\to B,\Theta)$ be a
generically smooth complex 
projective 
family of  stable pairs
 over a projective base variety ${B}$, defined over a 
field
$K$ of characteristic 0. Assume $({P}\to B,\Theta)$ is {\em  of maximal
 variation}, namely the morphism $B \to \obA_g$ is generically finite.  Let 
$G
\subset \Aut({P}\to B,\Theta)$ be a finite subgroup.    Then for a
sufficiently large integer 
$n$, the quotient pair
$(P_B^n/G,{\Theta}_{n}/G)$ is of logarithmic general
type. 
\end{proposition}

{\bf Proof.} Our first step is to replace $P \to B$ by a toroidal situation.
We have a moduli morphism $B \to \ofA_g$, hence at least a rational map to the
Voronoi compactification $B \das  \ofA'_g$, which is toroidal. Unfortunately we
might need finite covers to lift this to a morphism, since we are working with
stacks!   

In any case, there is a variety $B_1$ and  a  finite surjective  morphism $B_1
\to B$ such that $B_1 \to B \das \ofA'_g$ extends to a morphism. Taking the
normalization 
in 
the Galois closure of $K(B_1)$ over $K(B/G)$, we may assume that $B_1 \to B/G$
is Galois, with Galois group $G'$. The group 
$$G_1=G \mathop{\times}\limits_{\Aut B}G'$$
 acts on $P_1 = P \times_BB_1$ in such a way that $(P_1)^n_{B_1}/G_1=
P_B^n/G$. Thus we may replace $P \to B$ by $P_1 \to B_1$, in particular we may
assume that there is a morphism $B \to \ofA'_g$ lifting  $B \to \ofA_g$.

We can choose a $G$-equivariant resolution of singularities $B'\to B$ such that
$P' = P \times_BB'$ degenerates over a normal crossings divisor. Again, we
replace $P\to B$ by $P'\to B'$, so we may assume $B$ is nonsingular and $P \to
B$ degenerates over a normal crossings divisor. Denote by $U_B \subset B$ the
complement of this divisor, namely the smooth locus of $P \to B$, and
$U_P\subset P$ the inverse image. 

Since
$\ofA'_{g,1} \to  \ofA'_g$ is toroidal, we can apply \cite{A-Karu}, Lemma 6.2.
This implies that $(U_P\subset P) \to (U_B \subset B)$ is also a toroidal
morphism.  

Under this assumption we claim

\begin{claim} \begin{enumerate}
\item The pair $(P^n_B,\Theta_n)$ has log-canonical singularities, and 
\item The complement $P^n_B\setmin\Theta_n$  has only canonical singularities.
\end{enumerate}
\end{claim} 

{\bf Proof.} Observe that $({P}^n,{\Theta}_{n})$ is also
a family of  stable  pairs over a nonsingular base $B_1$,
satisfying 
 the toroidal assumptions.  Now the assertion follows directly from
\ref{Sssec:sqav-lc} and  \ref{Sssec:sqav-gorenstein}.\qed(Claim)

We go back to the proof of the proposition.

Since ${B}$ is generically finite over the moduli space and since the
generic fiber of $\pi$ is a smooth Abelian variety, Theorem 1.1 of
\cite{Kawamata:MMandKod} 
implies that
$$\det(\pi_*\omega_{{P}/{B}}^l) = \pi_*\omega_{{P}/{B}}^l$$ 
is a big invertible sheaf on
${B}$ for some positive integer
$l$.  Here we used the fact that 
$$h^0(\pi^{-1}(b), \omega_{\pi^{-1}(b)}^l) = 1 \quad
\forall b \in {B}.$$

Since
$\omega_{{P}/{B}}({\Theta})$ is relatively ample for $\pi$,
we conclude
$$\omega_{{P}/{B}}({\Theta}) \otimes \pi^*(\pi_*\omega_{{P}/{B}}^l)^n$$
is big for some sufficiently large integer $n$.  As we have the natural
inclusion map
$$\omega_{{P}/{B}}({\Theta}) \otimes \pi^*(\pi_*\omega_{{P}/{B}}^l)^n 
\hookrightarrow \{\omega_{{P}/{B}}({\Theta})\}^{1 + ln},$$
we conclude that the sheaf $\omega_{{P}/{B}}({\Theta})$ is
big.

Let $\Sigma
\subset {P}$ be the locus of fixed points of nontrivial elements of the group
$G$. 
Denote the sheaf of ideals
${\mathcal I}_{\Sigma}$.  For sufficiently large $n$ the sheaf
$$(\omega_{{P}/{B}}({\Theta}))^{\otimes n} \otimes \pi^*\omega_{ B} \otimes {\mathcal
I}_{\Sigma}^{|G|}$$
is big as well.  Taking the $n$-th fibered powers, we have that
$$(\omega_{{P}_{B}^n/{B}}({\Theta}_{n}))^{\otimes n} \otimes
\pi^*\omega_{B}^{\otimes n}
\otimes
\prod_{i = 1}^n p_i^{-1}{\mathcal I}_{\Sigma}^{|G|}$$
is a big sheaf on ${P}_{B}^n$.  Note that
$$\prod_{i = 1}^{n}p_i^{-1}{\mathcal I}_{\Sigma}^{|G|} \subset (\Sigma_{i =
1}^n p_i^{-1}{\mathcal I}_{\Sigma}^{|G|})^n,$$ 
where the latter ideal vanishes to order $\geq n|G|$ on the fixed points of
nontrivial elements of
$G$ in ${P}_{B}^n$.  Also note that
$$(\omega_{{P}_{B}^n/{B}}({\Theta}_{n}))^{\otimes n} \otimes
\pi^*\omega_{B}^{\otimes n} = 
(\omega_{{P}_{B}^n}({\Theta}_{n}))^{\otimes n}.$$
Therefore, for $l >> 0$ we have many invariant sections of
$\omega_{{P}_{B}^n}({\Theta}_{n})^{\otimes ln}$ vanishing on this fixed point 
locus to order $\geq ln \cdot |G|$. 
Let
$f:(V,D) \rightarrow ({P}_{B}^n,{\Theta}_{n})$  be an
equivariant good resolution of singularities with $f^{-1}({\Theta}_{n}) = D$.
 Now the following lemma, together with the Claim above, imply
that the invariant sections  of
$\omega_{{P}_{B}^n}({\Theta}_{n})^{\otimes ln}$ vanishing on the fixed point 
locus to order $\geq ln \cdot |G|$ descend to  sections of 
the pluri-log canonical divisors of a good resolution of the quotient pair
$({P}_{B}^n/G,{\Theta}_{n}/G)$, and hence it is of logarithmic general type.

%(Note that Claim 2.2.4 guarantees that a (invariant) 
%section of
%$(\omega_{{\mathfrak A}_{\mathfrak B}^n}(\mathfrak{\Theta}_{n}))^{\otimes ^k}$
%lifts to 
%that of
%$\omega_V(D)^{\otimes k}$.)

\begin{lemma}[\cite{A:stintegral}, lemma 4; see also \cite{C-H-M}, Lemma 4.1]
Let 
$V$ be a nonsingular  
variety with a normal 
crossing divisor $D$, and let $G \subset \Aut(V,D)$ be a finite subgroup of the
automorphism group of the  pair
$(V,D)$.  Let $\Sigma$ be the locus of the fixed points of the nontrivial
elements of $G$. 

Denote by $q:(V,D) \rightarrow (W = V/G,D_W = D/G)$ the natural morphism to
the quotient and choose  a good
resolution $r:({\tilde W},D_{\tilde W}) \rightarrow (W,D_W)$.

Then an invariant section
$$s \in H^0(V,\omega_V(D)^{\otimes k})^G$$
such that
$$s_x \in \omega_V(D)^{\otimes k} \otimes {\mathcal I}_{\Sigma}^{\otimes k|G|}
\hskip.1in
\forall x
\in
\Sigma \setmin D$$
comes from ${\tilde W}$, i.e., there exists
$$t \in H^0({\tilde W},\omega_{\tilde W}(D_{\tilde W})^{\otimes k})$$
such that $s = q^*r_*t$.
\end{lemma}

This completes the proof of Proposition \ref{Prop:quotients}, which implies
Theorem \ref{Th:lgt}, and Proposition \ref{Prop:correlation} follows.\qed

\subsection{A ``comprehensive'' family} In order to prove a result about {\em
all} principally polarized abelian varieties of a given dimension, we use the
following fact: there  exists 
a  projective family 
of principally polarized abelian varieties with  theta divisors
$\pi:({\fA}\rightarrow \fB,{\fXi})
$ over a quasi-projective base $\fB$ such that for any principally
polarized 
abelian variety with theta divisor
$(A,\Theta)$ defined over
$K$ there is a morphism $\Spec K \rightarrow \fB$ satisfying 

$$\begin{CD}
(A,\Theta) = \text{Spec} K \times (\fA,\fXi) @>>> (\fA,\fXi) \\
@VVV @VVV \\
\text{Spec} K @>>> \fB. \\
\end{CD}$$

Such a family can be found easily by noting that any principally
polarized abelian variety
$(A,\Theta)$ defined over $K$ can be embedded in a projective space using the
very ample linear 
system $|l\Theta|$ for some fixed $l \geq 3$, and therefore one can choose an
appropriate quasi projective subscheme of the Hilbert scheme  as the base
$\fB$.
Alexeev uses this fact during the construction of $\obA_g$, see
\cite{Alexeev:SQAV}.

\subsection{From correlation to uniformity}\label{Ssec:prolongable} Proposition
\ref{Prop:correlation} shows that the set $\cP$ of stably $S$-integral points
on $\fA\setmin \fXi$ is correlated. We now 
suggest an argument for uniformity,  replacing that  in Lemma
\ref{Lem:correlation}. We follow \cite{Pacelli:integral} and
\cite{A-Voloch:fibered}. The argument is complete in case $g\leq 2$.

Let $E_n \subset {\mathcal P}_\fB^n$ be the set of $n$-tuples of {\em distinct}
stably 
$S$-integral points.  Consider the tower of maps
$$\cdots \rightarrow E_{n+1} \rightarrow E_n \rightarrow E_{n-1}
\rightarrow \cdots \rightarrow E_1.$$
Let $E_n^{(m)}$ be the image of $E_m$ in $E_n$ and set $F_n^{(m)} =
\overline{E_n^{(m)}}$. We have a descending chain of closed subsets
	$$F_n^{(m)} \supset F_n^{(m+1)} \supset \cdots $$
which must stabilize to a closed set $F_n$, i.e., there exists $m_n$ such that
	$$F_n^{(m_n)} = F_n^{(m_n + 1)} = \cdots = F_n.$$
Accordingly, we have a tower of maps
$$ \cdots F_{n+1} \rightarrow F_n \rightarrow F_{n-1} \rightarrow \cdots
\rightarrow F_1.$$  
Ultimately we would like to conclude that all the $F_n$ are empty, which
implies 
the uniformity.

First observe that $F_{n+1} \rightarrow F_n$ is all surjective for all $n$.  In
fact, if we take $m \geq m_{n+1}, m_n$ then $E_{n+1}^{(m)} \rightarrow 
E_n^{(m)}$ is surjective by definition and hence $F_{n+1} =
\overline{E_{n+1}^{(m)}} \rightarrow F_n = \overline{E_n^{(m)}}$ is
surjective.  
Second, we wish to prove inductively that a fiber 
of $F_{n+1} \rightarrow F_n$ cannot have dimension $0, 1, ..., g$ and hence all
the $F_n$ are empty; in this paper we do this only for $g\leq 2$.  

  Denote by $\tau_{n+1}:F_{n+1} \rightarrow F_n$ be the
surjective morphisms as above.
It is enough to prove the following ``Inductive Statement $d$'' for all  $0
\leq d \leq g$:

\begin{quote}
{\bf Inductive Statement $\boldsymbol{d}$}: 
\em ``Suppose no fiber of $\tau_{n+1}:F_{n+1}
\rightarrow F_n$ 
has dimension 
$< d$ for all $ n \geq 1$, where $0 \leq d \leq g$.  Assume that the Lang-Vojta
conjecture holds true.  Then no fiber of 
$\tau_{n+1}:F_{n+1} \rightarrow F_n$ has dimension
$\leq d $ for all $n \geq 1$.''
\end{quote}

 Inductive Statement 0 can be proven without assuming any conjectures:

\begin{claim}
 No fiber of $\tau_{n+1}:F_{n+1} \rightarrow F_n$ has dimension 0 for all
$n \geq 1$.
\end{claim}

{\bf Proof of claim.}
Assume that a fiber of $\tau_{n+1}:F_{n+1} \rightarrow F_n$ has
dimension 0 for some $n$.  Then by the upper semicontinuity of the fiber
dimension, the dimension of a fiber must be 0 over some open subset $U$ of
$F_n$.  
Over $U$ we may assume the number of points in a fiber is also a constant $e$. 
The open set $U$ contains a point in $E_n^{(m)}$ where $m \geq m_n, m_{n+1}$
and $m - n > 
e$.  But as
$E_m$ consists of $n + (m - n)$-tuples of distinct points, this is
impossible.\qed(Claim)

We can prove  the ``upper end of the induction'' $d = g$ in general:

\begin{claim} Assume the Lang-Vojta conjecture holds true.
Suppose every fiber of $\tau_{n+1}:F_{n+1} \rightarrow F_n$ has dimension $g 
$ for all $n \geq 1$. Then $F_n = \emptyset$ for all $n$.
\end{claim}

{\bf Proof of claim.}
 Take an irreducible
component $M_n$ of $F_n$ and denote by $M_{n+k}$ its inverse images in
$F_{n+k}$ (which is also an irreducible component in $F_{n+k}$) for $k \geq
1$.  Then $M_{n+k} = {M_{n+1}}_{M_n}^k$, where $M_{n+1} \rightarrow M_n$ is a
family of principally polarized quasi-abelian varieties of dimension $g$, with 
theta divisors
${\Theta}_{M_{n+1}}$.  By Proposition \ref{Prop:correlation}, for sufficiently
 large 
 integer  $k$, the set of $k$-tuples of  stably
$S$-integral points ${\mathcal P}_{M_n}^k$ is not dense in the fibered power
$({M_{n+1}})_{M_n}^k$.  On the  other hand,
$M_{n+k}$ has {\em by definition} a dense set of  stably
$S$-integral points, a contradiction.\qed(Claim)

From now on fix $d\geq 1$, We now make a few general reduction steps for
Inductive Statement $d$. We note 
that in general we have  $F_{n+k}
\subset {F_{n+1}}_{F_n}^k$ for $k \geq 1$. Consider an irreducible component
$B\subset  F_n$ over which $X = \tau_{n+1}^{-1}B$ has relative dimension $d$. 
 Consider $\tau =
\tau_{n+1}|_X : X \to B$, and 
$\cQ = E_{n+1}\cap X$. In terms of Definition \ref{Def:strong-correlation}, we
can reformulate our problem as follows:  

\begin{lemma}
In order to prove Inductive Statement $d$, it suffices to show that $\cQ$ is
strongly $k$-correlated for some integer $k$.
\end{lemma}

Indeed, we have that the  dimension of every fiber of $F_{n+k} \to F_n$ is
$\geq kd$.\qed

Denote by $X_i, i=1,\ldots,c'$  the irreducible components of $X$ of relative
dimension $d$ over $B$. Let $X_i'\to X_i$ be the normalization, and $X_i' \to
B_i'\to B$ the Stein factorization. In order to keep things inside a family of
smooth quasi-abelian varieties, consider $A = B' \times_\fB\fA$ and denote by
$X_i''$ the image of $X_i'$ in $A$. Note that we have a factorization $X_i' \to
X_i''\to X_i$. Consider the subset $\cQ_i'$ as in Proposition
\ref{Prop:strong-correlation}. Writing $\Theta = B'\times_\fB\fXi$, we can view
$\cQ_i'$ as a subset of the set of 
stably integral points on $A\setmin \Theta$. 

By Proposition \ref{Prop:strong-correlation}, we have the following:

\begin{lemma}\label{Lem:sub-family} In order to prove Inductive Statement $d$,
it suffices to show 
that $\cQ_i'\subset X_i''(K)$ is correlated for $ X_i''(K)\to B'$ for all $i$. 
\end{lemma}

What can we say about correlation of $\cQ_i'$? 
For simplicity of notation denote $X_i'' = H \subset A$, and replace $B$ by
$B'$. We may even replace $B$ by a nonempty 
open subset, therefore we may assume $H \to B$ is flat. We can now classify the
family $H \to B$. We apply the
following proposition, which is essentially due to Ueno \cite{Ueno:classI} and
Kawamata \cite{Kawamata:Bloch}, to
our situation. 

\begin{proposition} Let ${H} \rightarrow {B}$ be a flat
family of geometrically irreducible closed subvarieties of a family of
smooth quasi-abelian varieties 
${A} 
\rightarrow {B}$,
defined over a number 
field $K$,  where $A$  is a torsor on an abelian scheme $\tilde{A}\to
B$.  Then there exists an
abelian subscheme
${D}
\rightarrow {B}$ such that ${H}$ is invariant under
translation 
by ${D}$ and such that the quotient
${H}/{D} \rightarrow {B}$ is a family of
varieties of 
general type.
\end{proposition}

 Ueno and Kawamata state the above proposition only  in the case ${B} = 
\Spec \CC$. The general case is given in \cite{McQuillan:complements}.

Consider the family ${H}/{D} \rightarrow {B}$. It is a family of
varieties of general type. We have two cases to consider

Case 1: ${H}/{D} \rightarrow {B}$  has relative dimension $> 0$.

By the Fibered Power Theorem (Theorem \ref{Th:fibered-power}), $(H/{D})^k_{B}$
dominates a positive 
dimensional  variety of general type, for some integer $k$. This immediately
implies that $H^k_{B}$ dominates a positive 
dimensional  variety of general type. Lang's conjecture implies that
$H^k_{B}(K)$ is not dense, which implies that $\cQ_i'$ is
$k$-correlated, which is what we wanted.

Case 2: ${H}/{D} \rightarrow {B}$ is an isomorphism.

In this case, the generic fiber of $H \subset {A}$ is a 
 translate of the generic fiber of ${D}$, which is an abelian variety.
 {\em The issue here is whether or not the set of stably
$S$-integral points contained in this sub-family is correlated.}

\subsection{The case of an elliptic subscheme}
In this section we provide an argument for case $d=1$, completing the
proof of Inductive Statement $1$. The argument uses
Pacelli's strong uniformity  
results for elliptic curves 
(see \cite{Pacelli:integral}). This completes the proof of the Main Theorem. 
At the end of the paper we  
 discuss a possible  line of argument which could lead to a proof in
arbitrary dimension.

Thus we assume $\tau:H \rightarrow B$ is a family of elliptic
curves inside of principally polarized abelian varieties. For simplicity of
notation we denote $\cQ_i' = \cR$. 

Take a point $P \in\cR$, and let $b =\tau(P) \in B$. 
Choose a point $O \in E \cap \Theta$ (which may not be $K$-rational).  Note
that 
there is a constant
$a$ which only depends on the family $H \rightarrow B$ and
 such that $a \geq \length_K
{\mathcal O}_{E \cap 
\Theta}$ and hence that the extension degree of the residue field at $O$ over
$K$ is bounded by this number $[k(O):K] \leq a$.  Let $(E \times \Spec k(O),
O)$ be 
the elliptic 
curve extended by extension to the residue field $k(O)$.  Then the following
lemma implies 
that $P \in (E \times \Spec k(O) \setmin 
\{O\})(k(O))$ is a stably $S_{k(O)}$-integral point, where $S_{k(O)}$ is the
set of places in $k(O)$ over $S$. 

\begin{lemma} Let $(A,\Theta)$ be aprincipally polarized abelian variety and 
$E \subset A$ be an  
elliptic curve with the origin $O \in E \cap \Theta \subset A$ defined over a
number field 
$L$. Let $Q \in E \setmin \Theta \subset
A \setmin \Theta$ be 
an $L$-rational point, i.e., $Q \in (E \setmin 0)(L) \subset (A \setmin
\Theta)(L)$, 
which is a stably 
$S_L$-integral point of the pair $(A,\Theta)$.

Then $Q$ is also a stably $S_L$-integral point of the pair $(E,O)$.
\end{lemma}

{\bf Proof.} First remark that by Proposition \ref{Prop:stint-any-L} we
are allowed to 
take any finite extension of $L$ in order to prove the assertion.

By taking such a  finite extension of $L$, we may assume that 
$(A\rightarrow \Spec L,\Theta)  
$ has the  stable quasi-abelian model 
$\pi:({\mathcal P}\to \Spec
{\mathcal O}_{L,S_L},\tilde{\varTheta})$, and that $(E,O) \rightarrow \Spec L$
has a   stable model ${\mathcal E} \rightarrow \Spec {\mathcal
O}_{L,S_L}$. 

 For any prime $\fp \in \Spec {\mathcal
O}_{L,S_L}$, we need to show that the point $\overline{Q}_{\mathcal E}$ is
disjoint from the origin 
$O$ in the 
fiber
${\mathcal E}_\fp$ over $\fp$, where $\overline{Q}_{\mathcal E}$ is the Zariski
closure of the $L$-rational 
point $Q$ in ${\mathcal E}$.  Therefore, instead of working over
${\cO}_{L,S_L}$, we may work over the completion $R$ of ${\cO}_{L,S_L}$ at
the prime $\fp$. 

We denote by $\eta$ the generic point of $\Spec R$. In a slight abuse of
notation, we use the same letters for the objects over $\Spec R$ as for those
over $\Spec {\mathcal O}_{L,S_L}$.      

Replacing $R$ by a finite extension we may assume there exists a section $s:
\Spec R
\rightarrow \cP$ so that its image $O_s$ sits
in the locus where 
the morphism $\pi$ is smooth. This can be done by picking a smooth point on the
central fiber $\cP_\fp$ and lifting using Hensel's lemma.

Let ${\mathcal N}_A \rightarrow \Spec R$ be the N\'eron model of $A$
(with the origin $O_s$), and ${\mathcal N}_E \rightarrow \Spec R$ the N\'eron model
of $E$ (with 
the origin $O$).

By the universal property of the N\'eron model, we have a morphism
$$\phi:{\mathcal N}_E \rightarrow {\mathcal N}_A$$ extending the inclusion $E 
\hookrightarrow A$.
By Proposition \ref{Prop:sqav-neron} we also have a morphism 
$$\psi:{\mathcal N}_A \rightarrow {\cP},$$
extending the isomorphism
$$({\mathcal N}_A)_{\eta} \overset{\sim}{\rightarrow} {\mathcal P}_{\eta} =
A.$$ 
Therefore, we have a morphism
$$\psi \circ \phi:{\mathcal N}_E \rightarrow {\mathcal P}$$
extending the inclusion
$$E \hookrightarrow A.$$

Under this morphism, the origin $O_{{\mathcal N}_E}$ maps into
$\tilde{\varTheta}$, while 
$\overline{Q}_{{\mathcal 
N}_E}$ maps to $\overline{Q}_{{\mathcal P}}$,

Denote by $\overline{Q}_{{\mathcal
N}_E}$ and $\overline{Q}_{{\mathcal P}}$ are the closures of the $L$-rational
point $Q$ in 
${\mathcal N}_E$ and ${\mathcal P}$, respectively. Then clearly $\psi \circ
\phi(\overline{Q}_{{\mathcal  
N}_E}) = \overline{Q}_{{\mathcal P}}$.

Since $Q$ is a
stably $S_L$-integral point of the pair $(A,\Theta)$, we have that
$\tilde{\varTheta}$ and $\overline{Q}_{{\mathcal P}}$ are disjoint in
${\mathcal P}$. A-fortiori,  $O_{{\mathcal N}_E}$ and $\overline{Q}_{{\mathcal 
N}_E}$ are disjoint in ${\mathcal N}_E$. 

Let ${\mathcal E}$ be  the stable model of $E$ over  $\Spec R$.
Note that ${\mathcal N}_E$ and ${\mathcal E}$
 are isomorphic in a neighborhood of $O_{{\mathcal N}_E}$. We conclude that
$O_{{\mathcal E}}$ 
and 
$\overline{Q}_{{\mathcal E}}$ are disjoint, and hence that $Q$ is a stably
$S_L$-integral point of the pair $(E,O)$.
\qed(Lemma)

Going back to Inductive Statement 1, recall that the extension degree $[k(O):K]
\leq a$ over the fixed 
number field $K$ is uniformly bounded. In particular, $[k(O):\QQ]\leq d =
a\cdot 
[K:\QQ]$.  Pacelli's result \cite{Pacelli:integral} asserts that, 
assuming the Lang-Vojta conjecture, there is a uniform bound $N(d,S)$ such that
for 
any elliptic curve $(E,\{O\})$ defined over a number field $L$ of 
degree $\leq d$, 
 the number of the stably $S_L$-integral points is
uniformly bounded 
$$\# E(L,S_L)^{stable} < N(d,S).$$

In our situation, this implies that the number of stably $S$-integral points on
$A\setmin \Theta$ lying in $E$ is uniformly bounded by $N(d,S)$, which in
particular says 
that the points are $N(d,S)$-correlated.

This completes the proof of  Inductive Statement $1$, which completes the 
case $\dim A = g = 2$ and hence the Main Theorem.
  \qed

\subsection{Towards higher dimensions} Finally we discuss a possible line of
argument for higher dimensions (which would also lead to a result about
arbitrary polarizations).  

If one considers Lemma \ref{Lem:sub-family} and Case 2 in the discussion
following that Lemma, one reduces to the following conjecture: 

\begin{conjecture} Let $(A\to B, \Theta)$ be a  family of smooth principally 
polarized quasi-abelian varieties over a number field $K$. Let $H\subset A$ be
a family of quasi-abelian subvarieties. Let $\cP\subset A(K)$ be the set of
stably integral points. Then $\cP\cap H$ is correlated with respect to $H \to
B$.
\end{conjecture}

We would like, at least, to show that this conjecture follows from the
Lang-Vojta conjecture. 

Denote  $\Theta_H = H \cap \Theta$. Replacing $B$ by a nonempty open subset, we
may assume $\Theta$ does not contain a fiber of $H \to B$. 
Then $(H\to B, \Theta_H)$ can be viewed as a  family of polarized quasi-abelian
varieties. The issue is, that these are not necessarily principally polarized.

In the recent preprint \cite{Alexeev:toric}, Alexeev defines a complete moduli
space for such pairs as well. We call these ``Alexeev stable pairs''
below. This suggests the following approach to the 
problem: 

\begin{enumerate}
\item Define ``Alexeev stably integral points'' of $(H,\Theta_H)$ to be
rational points
which are integral on the complement of $\Theta_H$ in an Alexeev stable model
of a pair $(H, \Theta_H)$. 
\item Give a criterion for Alexeev stably integral points in terms of N\'eron
models. 
\item Deduce that a stably integral point of a pair $(A,\Theta)$ is also stably
integral on $(H,\Theta_H)$.
\item Assuming Lang-Vojta, reduce the problem to a problem on moduli of $n$-pointed Alexeev stable
pairs similar to Theorem \ref{Th:lgt} and Proposition \ref{Prop:quotients}.
\item Prove a result analogous to Proposition \ref{Prop:quotients}.
\end{enumerate}

All but the last step seem straightforward. The main issues in the last step
are:
\begin{enumerate} 
\item Suppose $(P\to B, \Theta)$ is  an Alexeev stable pair of maximal
variation, defined over a field $K$, over a projective irreducible  nonsingular
base $B$, with smooth generic 
fiber. There exists $\boldsymbol{\varepsilon}>0$ such that, for all $n$, the
pair $(P^n_B, 
\boldsymbol{\varepsilon}\Theta_n)$ has log-canonical singularities. 
\item For such  $(P\to B, \Theta)$, the sheaf $\omega_{P/B}^m(\Theta)$ is big
for some $m>0$.
\end{enumerate}

We expect that these statements can be proven using Alexeev's work.

%\bibliographystyle{alpha}
%\bibliography{cite}
%
%\end{document}

%
%------------------------------------------------------------

\end{document}